\documentclass[a4paper,12pt]{article}
\usepackage{amsmath,amssymb,graphicx}
\usepackage{algorithm}
\usepackage{algorithmic}
\usepackage{caption,color}
\usepackage{hyperref} 
\usepackage{tikz}
\usepackage{subcaption}

\title{Spectral extremal graphs for fan graphs\thanks{Email: 
\url{Ylj99@hnu.edu.cn} (L. Yu), 
\url{ytli0921@hnu.edu.cn} (Y. Li), \url{ypeng1@hnu.edu.cn} (Y. Peng, corresponding author)}}
\author{\small{Loujun Yu$^1$, Yongtao Li$^2$, Yuejian Peng$^1$}\\
{\small $^1$School of Mathematics, Hunan University}\\
{\small Changsha, Hunan, 410082, P.R. China} \\ 
{\small $^2$School of Mathematics and Statistics, Central South University} \\ 
 {\small Changsha, Hunan, 410083, P.R. China}  
\\\makeatletter
}
\newcommand\tabcaption{\def\@captype{table}\caption}
\date{\today} 
\makeatother

\newtheorem{theorem}{Theorem}[section]
\newtheorem{thm}[theorem]{Theorem}

\newtheorem{lemma}[theorem]{Lemma}

\newtheorem{conjecture}[theorem]{Conjecture}
\newtheorem{conj}[theorem]{Conjecture}
\newtheorem{claim}{Claim}

\newcommand*{\QED}{\hfill\ensuremath{\square}}

\begin{document}
\maketitle
\begin{abstract}
A well-known result of Nosal states that 
 a graph $G$ with $m$ edges and 
$\lambda(G) > \sqrt{m}$ contains a triangle. 
Nikiforov [Combin. Probab. Comput. 11 (2002)] extended this result to cliques 
by showing that if $\lambda (G) > 
\sqrt{2m(1-1/r)}$, then $G$ contains a copy of $K_{r+1}$. 
Let $C_k^+$ be the graph obtained from a cycle $C_k$ 
by adding an edge to two vertices with distance two, and let $F_k$ be the friendship graph consisting of 
$k$ triangles that share a common vertex.  
Recently, Zhai, Lin and Shu [European J. Combin. 95 (2021)], Sun, Li and Wei [Discrete Math. 346 (2023)], and 
Li, Lu and Peng [Discrete Math. 346 (2023)] 
proved that if $ m\ge 8$ 
and $\lambda (G) \ge \frac{1}{2}
(1+\sqrt{4m-3})$, then $G$ contains a copy of 
$C_5,C_5^+$ and $F_2$, respectively, unless {$G=K_2\vee 
\frac{m-1}{2}K_1$}. 
In this paper, we give a unified extension 
by showing that such a graph 
contains a copy of $V_5$, where 
$V_5=K_1\vee P_4$ 
is the join of a vertex and a path on four vertices. 
Our result extends the aforementioned results 
since $C_5,C_5^+$ and $F_2$ 
are proper subgraphs of $V_5$. 
In addition,  
we prove that if $m\ge 33$ and $\lambda (G) 
\ge 1+ \sqrt{m-2}$, then $G$ contains a copy of $F_3$, unless $G=K_3\vee \frac{m-3}{3}K_1$. 
This confirms a conjecture 
on the friendship graph  $F_k$  
in the case $k=3$. Finally, we conclude some spectral extremal graph problems concerning the large fan graphs and wheel graphs. 
\end{abstract}



\section{Introduction}

Let $K_n$ be the complete graph on $n$ vertices, and
$K_{s,t}$ be the complete bipartite graph with parts of sizes
$s$ and $t$.
We write  $C_n$ and $P_n$ for the cycle and
path on $n$ vertices respectively. Let $kG$ be the union 
of $k$ disjoint copies of $G$. 
Let $G \vee H$ be the join graph obtained from $G$ and $H$ by joining  each vertex of $G$
 to each vertex of $H$.  

Dating back to 1970,  Nosal \cite{Nosal1970} (see, e.g., \cite{Ning2017-ars} for alternative proofs) showed that for 
every triangle-free graph $G$ with $m$ edges, we have  
$\lambda (G)\le \sqrt{m} $. 
Later, Nikiforov \cite{Niki2002cpc,Niki2006-walks,Niki2009jctb} extended Nosal's result by 
proving that 
if $G$ is $K_{r+1}$-free, then 
\begin{equation}       \label{eq1}
\lambda (G) \le \sqrt{2m\Bigl( 1-\frac{1}{r}\Bigr)}. 
\end{equation}
Moreover, the equality holds 
if and only if 
$G$ is a complete bipartite graph for $r=2$, 
or a regular complete $r$-partite graph for each $r\ge 3$. 
Tur\'{a}n's theorem 
can be derived from  (\ref{eq1}).  
Indeed, using Rayleigh's inequality, we have 
$\frac{2m}{n}\le \lambda (G)\le  (1-\frac{1}{r})^{1/2} \sqrt{2m}$, 
which yields $ m \le (1- \frac{1}{r}) \frac{n^2}{2}$.  
Thus,  (\ref{eq1}) could be viewed as  
a spectral extension of Tur\'{a}n's theorem. 
Moreover,  since each tree is triangle-free, 
(\ref{eq1}) implies a result of Lov\'{a}sz and Pelik\'{a}n \cite{LP1973}, 
which asserts that if $G$ is a tree on $n$ vertices, then $\lambda (G)\le \sqrt{n-1}$, 
with equality  if and only if $G$ is a star.  
Furthermore, Bollob\'{a}s and Nikiforov \cite{BN2007jctb}  
 conjectured that  
if  $G$ is a $K_{r+1}$-free graph 
of order at least $r+1$ 
with $m$ edges, then 
\begin{equation} \label{eq-BN}
{  \lambda_1^2(G)+ \lambda_2^2(G) 
\le 2m\Bigl( 1-\frac{1}{r}\Bigr) }. 
\end{equation} 
This conjecture was firstly confirmed
by Ando and Lin \cite{AL2015} for $r$-partite graphs.  
In 2021,  the case $r=2$ was confirmed by Lin, Ning and Wu \cite{LNW2021} 
by using doubly stochastic matrix theory. 
Later,  Nikiforov  \cite{Niki2021} provided a simple proof by using elementary inequalities.  
In 2022, Li, Sun and Yu \cite{LSY2022} proved 
an extension for $\{C_3,C_5,\dots ,C_{2k+1}\}$-free graphs.  
 In 2024, Zhang \cite{Zhang2024} proved the conjecture for all regular graphs. 
 Nevertheless this intriguing 
 problem remains open, 
we refer the readers to 
\cite{ELW2024} 
for two related conjectures.  
As an application of the triangle case of (\ref{eq-BN}), 
Lin, Ning and Wu \cite{LNW2021} stimulated 
an interesting question to investigate the spectral 
problems for non-bipartite triangle-free graphs; 
see \cite{ZS2022dm,LLH2022, 
LP2022second,LFP2023-solution,Wang2022DM} 
for some recent developments. 

Both (\ref{eq1}) and (\ref{eq-BN}) 
boosted the great interest of studying the maximum spectral radius
for an $F$-free graph with given number of edges. 
For example, see \cite{Niki2009laa,ZS2022dm,Wang2022DM} for  $C_4$-free graphs,  
\cite{ZLS2021} for  $K_{2,r+1}$-free graphs,  
 \cite{ZLS2021,MLH2022} for $C_5$-free or $C_6$-free graphs, 
\cite{LLL2022} for $C_7$-free graphs, 
\cite{LSW2022,FYH2022,LW2022} for $C_5^{+}$-free and $C_6^{+}$-free graphs and in general,  \cite{LZS2024} for $C_{k}^+$-free graphs,  
where $C_k^{+}$ 
is a graph on $k$ vertices obtained from $C_k$ 
by adding a chord between two vertices with distance two; 
see \cite{Niki2021} for  $B_k$-free graphs, 
where $B_k$ denotes the book graph which consists  of $k$ triangles 
sharing a common edge, 
\cite{LLP2022} for $F_2$-free graphs with given number of edges, 
where $F_2$ 
is the friendship graph  consisting of two triangles  intersecting in 
a common vertex. 
In particular, Zhai, Lin and Shu \cite{ZLS2021} 
proved the following result for $C_5$-free or $C_6$-free graphs. 

\begin{thm}[Zhai--Lin--Shu \cite{ZLS2021}, 2021] 
\label{thm-ZLS}
If $G$ is a $C_5$-free graph with $m\ge 8$ edges, then 
\[  \lambda (G) \le \frac{1+\sqrt{4m-3}}{2}, \]
where the equality holds if and only if $G=K_2 \vee \frac{m-1}{2}K_1 $. 
\end{thm}

Recall that $C_t^{+}$
{denotes} the graph on $t$ vertices obtained from 
$C_t$ by adding a chord between two vertices of distance two.   
Observe that $C_4^{+} = B_2$ is the book graph on $4$ vertices, 
and it was studied in \cite{Niki2021}. 
Note that $C_5^{+} \nsubseteq K_2 \vee \frac{m-1}{2}K_1 $, 
but $C_5^{+} \subseteq K_3 \vee \frac{m-3}{3}K_1$. 
The following result was proved by Li, Sun and Wei \cite{LSW2022}.

\begin{thm}[Sun--Li--Wei \cite{LSW2022}, 2023] 
\label{thm-LSW2022}
If $G$ is a $C_5^{+}$-free graph with $m\ge 8$ edges, then  
\[  \lambda (G) \le \frac{1+\sqrt{4m-3}}{2},  \]
equality holds if and only if $G=K_2 \vee \frac{m-1}{2}K_1 $. 
\end{thm}

Observe that $C_t\subseteq C_t^{+}$, so Theorem \ref{thm-LSW2022} 
can imply Theorem \ref{thm-ZLS}.    

\begin{theorem}[Li--Lu--Peng \cite{LLP2022}, 2023]  \label{thm-m-F2-36}
If $m\ge 8$ and $G$ is an $F_2$-free graph with $m$ edges, then 
\[  \lambda (G)\le \frac{1+\sqrt{4m-3}}{2}, \]
where the equality holds if and only if $G=K_2\vee \frac{m-1}{2}K_1 $. 
\end{theorem}

Let $V_k=K_1\vee P_{k-1}$ be the fan graph on 
$k$ vertices. 
It is worth noting that 
$V_3$ is a triangle and $V_4$ is a book on $4$ vertices, which were studied by Nosal \cite{Nosal1970}, Nikiforov \cite{Niki2002cpc,Niki2021} and 
Ning and Zhai \cite{NZ2021}. 
In addition, we can see that $V_5$ 
is the power of a path $P_5$. 
Recently, the $n$-vertex $V_5$-free graph 
with maximal spectral radius 
was also studied by 
Zhao and Park \cite{ZP2022}. 
Observe that every $C_5$-free, $C_5^+$-free or 
$F_2$-free graph must be $V_5$-free.  
In this paper, we shall investigate 
the $V_5$-free graphs with given size 
$m$ (instead of the order $n$), and we 
give a unified extension on 
Theorems \ref{thm-ZLS}, \ref{thm-LSW2022} and 
\ref{thm-m-F2-36}.

\begin{theorem}\label{Thm-V5}
If $G$ is a $V_5$-free graph with $m\ge  8$ edges, then 
\[ \lambda (G)\le \frac{1+\sqrt{4m-3}}{2}, \] 
 and the equality holds if and only if $G = K_2\vee \frac{m-1}{2} K_{1}$. 
\end{theorem}

Recently, Li, Lu and Peng \cite{LLP2022} 
 proposed a conjecture for $F_k$-free graphs.  

\begin{conjecture}[Li--Lu--Peng, 2023]  
 \label{conj-edges-Fk} 
Let $k\ge 2$ be fixed and $m$ be large enough. 
 If $G$ is an $F_k$-free graph with $m$ edges, then 
\[  \lambda (G)\le \frac{k-1 +\sqrt{4m -k^2+1}}{2}, \]
and equality holds if and only if 
$G=K_k \vee\frac{1}{k}\left(m-\binom{k}{2} \right) K_1$. 
\end{conjecture}

In this paper, we confirm Conjecture \ref{conj-edges-Fk}
 in the case $k=3$. 

\begin{theorem}\label{Thm-F3}
If $G$ is an $F_3$-free graph with $m\ge  33$ edges, then 
\[ \lambda (G)\le 1+\sqrt{m-2}, \] 
 and the equality holds if and only if $G = K_3\vee \frac{m-3}{3} K_1$.
\end{theorem}

\section{Preliminaries}
 Throughout this paper, we use the following {notations}. 
{For each forbidden graph $F$, let $G^*$ denote an extremal graph attaining maximum spectral radius among all $F$-free graphs with $m$ edges. An extremal vertex of $G^*$, denoted by $u^*$, is a vertex corresponding {to} the maximum coordinate of the Perron vector of $G^*$.}
 Given a vertex $u\in V(G)$, 
 {let $N_G(u)$ be the set of neighbors of $u$,
and $d_G(u)$ be the degree of $u$ in $G$.} 
Moreover, we denote  $N_G[u]=N_G(u)\cup\{u\}$. 
Sometimes, we will eliminate the subscript and write {$N(u)$ and $d(u)$} if not necessary. 
For a subset $U\subseteq V(G)$, we write $e(U)$ for the
number of edges with two endpoints in $U$. 
For two disjoint sets $U,W$, we write $e(U,W)$ for the number of edges between $U$ and $W$. 
{For simplicity,  we denote by 
$N_{U}(u)$ the set of vertices of $U$ that are  adjacent to $u$, i.e, $N_U(u)=N_G(u)\cap U$, and let $d_U(u)$ be the number of vertices of $N_U(u)$.}

\begin{lemma}[Wu--Xiao--Hong \cite{WXH2005}, 2005] \label{le-niki}
Let $G$ be a connected graph
and $(x_1,\ldots ,x_n)^T$ be a Perron vector of $G$,
where the coordinate $x_i$ corresponds to the vertex $v_i$.
Assume that
 $v_i,v_j \in V(G)$ are vertices such that $x_i \ge x_j$, and $S\subseteq N_G(v_j) \setminus N_G(v_i)$ is non-empty.
 Denote $G'=G- \{v_jv : v\in S\} +
\{v_iv : v\in S\}$. Then $\lambda (G) < \lambda (G')$.
\end{lemma}

\begin{lemma}[Zhai--Lin--Shu \cite{ZLS2021}, 2021]\label{le-zhai}
If $F$ is a $2$-connected graph, then
$G^*$ is connected. Moreover, there is no cut vertex in $V(G^*)\setminus\{u^*\}$. Furthermore, 
we have $d(u)\ge  2$ for any vertex $u\in V(G^*)\setminus N[u^*]$. 
\end{lemma}

Let $R_{s,t}$ be the graph obtained from 
$s$ copies of $K_4$ and $t$ independent edges sharing a common vertex. Namely, we have $R_{s,t}=K_1 \vee 
(sK_3 \cup tK_1)$. 

\begin{lemma}[Zhai--Lin--Shu \cite{ZLS2021}, 2021]\label{le-zhai-1}
If $k\ge 1$ and $m=6s+t$, then 
\[ 
\lambda(R_{s,t})< \frac{1+\sqrt{4m-3}}{2}. \] 
\end{lemma}

\begin{lemma}[Liu--Peng--Zhao \cite{Liu2006}, 2006]\label{le-Liu}
A graph $G$ contains no induced subgraph isomorphic to $P_6$ if and only if each connected induced subgraph of $G$ contains a dominating induced $C_6$ or a dominating complete bipartite graph.
\end{lemma}

\begin{lemma}[Erd\H{o}s--Gallai \cite{Erdos1959}, 1959]\label{le-E-G}
 Let $G$ be a $(k+1)K_2$-free graph of order $n$. Then 
 \[ e(G)\le \max \left\{ \binom{2k+1}{2},
 {\binom{k}{2}} + (n-k)k \right\}. \]
 and the equality holds if and only if $G = K_{2k+1}$ or $K_{k}\vee (n-k)K_1$.
\end{lemma}

\section{Proof of Theorem \ref{Thm-V5}}

In this section, we will consider the maximum spectral radius of $V_5$-free graphs with $m$ edges. 
{Recall that  $G^*$ is a spectral extremal graph for the fan graph $V_5$, and $u^*$ is an extremal vertex of $G^*$.} For convenience, we will use {$\lambda$} and $\mathbf{x}$ to denote the spectral radius and the Perron vector of $G^*$. 
Let $U=N(u^*)$ and $W=V(G^*)\setminus N[u^*]$.  {Since $V_5$ is 2-connected, by Lemma \ref{le-zhai}, we know that $G^*$ is connected. In addition, }
we get $\lambda \ge  \lambda (K_2\vee \frac{m-1}{2} K_1)=\frac{1+\sqrt{4m-3}}{2}$ since $K_2\vee \frac{k-1}{2}K_1 $ is $V_5$-free.
 Hence $\lambda^2-\lambda  \ge  m-1$. 
 Since 
 \[ \lambda^2 x_{u^*} = |U| x_{u^*} +
\sum_{u\in U} d_U(u)x_u + \sum_{w\in W} d_U(w) x_w, \] 
which together with $\lambda x_{u^*}=\sum_{u\in U}x_u$  yields
\[  (\lambda^2 - \lambda ) x_{u^*}
= |U| x_{u^*} + \sum_{u\in U} (d_U(u) -1)x_u
 + \sum_{w\in W} d_U(w) x_w.  \] 
Note that $\lambda^2 - \lambda \ge m-1
= |U| + e(U) + e(U,W) +e(W) -1$. Then
  \begin{align}\label{eqc}
   (e(U)+e(W)+e(U,W)-1)x_{u^*}\le  \sum\limits_{u\in U}(d_U(u)-1)x_u+\sum\limits_{w\in W}d_U(w)x_w. 
  \end{align}
By simplifying, we can get 
\begin{equation} \label{eqa}
 e(W)\le \sum\limits_{u\in U}(d_U(u)-1)\frac{x_u}{x_{u^*}} 
 -  e(U) + \sum_{w\in W} d_U(w) \frac{x_w}{x_{u^*}} - e(U,W) +1, 
 \end{equation}
where the equality holds if and only if $\lambda^2- \lambda =m-1$.

\begin{claim} \label{claim-eU}
$e(U)\ge 1$.
\end{claim}

\noindent \textbf{Proof}. 
As $\lambda \ge \frac{1+\sqrt{4m-3}}{2}> \sqrt{m}$ when $m>1$, we have
\begin{align*}
mx_{u^*}& < \lambda^{2}x_{u^*}=|U|x_{u^*}+\sum\limits_{u\in U}d_U(u)x_u+\sum\limits_{w\in W}d_U(w)x_w \\ 
&\le |U| x_{u^*} + 2e(U) x_{u^*} + e(U,W)x_{u^*}.
\end{align*}
Observe that $m=|U| + e(U) + e(U,W) + e(W)$. 
Therefore, the above inequality  implies that $e(W)<e(U)$. If $e(U)=0$, then $e(W)<0$, a contradiction. Thus, $e(U)\ge  1$. \QED

\medskip 

An isolated vertex in $G^*[U]$ is called a trivial component. 
Let $U_0$ denote the set of isolated vertices in $G^*[U]$. 
By Claim \ref{claim-eU},  
there {exists} at least one non-trivial  {component} in $G^*[U]$.
Let $\mathcal{H}$ be the set of all non-trivial components in $G^*[U]$.  For each non-trivial component $H\in \mathcal{H}$, we define 
\begin{equation} \label{eq-eta}
 \eta_1(H) :=\sum\limits_{u\in V(H)}(d_H(u)-1)\frac{x_u}{x_{u^*}}-e(H). 
 \end{equation} 
 Clearly,  we have 
$\eta_1 (H) \le e(H) - |V(H)| $.  
Then (\ref{eqa}) implies 
\begin{align}\label{eqb}
 e(W) \le \sum\limits_{u\in U}(d_U(u)-1)\frac{x_u}{x_{u^*}} 
 -  e(U) +1 
 = \sum\limits_{H\in \mathcal{H}}\eta_1(H)-\sum\limits_{u\in U_0}\frac{x_u}{x_{u^*}}+1,  
\end{align}
with equality if and only if 
$\lambda^2- \lambda =m-1$ and $x_w=x_{u^*}$ 
for any $w\in W$  with $d_U(w)\ge 1$. 

 Recall that $G^*$ is $V_5$-free, we  know that $G^*[U]$ is $P_4$-free. 
 It implies that any non-trivial component of $\mathcal{H}$ is isomorphic to a triangle or a star. Observe that 
 $\eta_1 (K_3)\le 0$ and $\eta_1 (K_{1,s}) \le -1$ 
 for each $s\ge 1$. 
 Next, we will prove that $e(W)=0$, and then we show that  not all components of $G^*[U]$ are isomorphic to $K_3$.

 \begin{claim} \label{eq-W0}
     If $W\neq \varnothing$, then 
     $e(W)=0$. 
 \end{claim} 

\noindent \textbf{Proof}. 
First of all,  
for any component $H\in  \mathcal{H}$, 
we know that $H$ is isomorphic to 
$K_3$ or $K_{1,s}$ for some $s\ge 1$. Then 
\[ \eta_1(H)=\sum\limits_{u\in V(H)}(d_H(u)-1)\frac{x_u}{x_{u^*}}-e(H) 
\le  0, \] 
 and the equality holds if and only if {$H= K_3$ and } $x_u=x_{u^*}$ for any $u\in V(H)$. 
 By (\ref{eqb}), we get $e(W)\le 1$. 
 Suppose on the contrary that $e(W)\ge  1$, then  $e(W)=1$. Consequently, we have
 $\eta_1 (H)=0$ for any $H\in \mathcal{H}$, and so 
 $H$ is isomorphic to $K_3$. 
 Moreover, the equality in (\ref{eqb}) leads to 
 $U_0=\varnothing $ and {$G^*[U]= t K_3$} for some $t\ge 1$.  
In addition, we have $x_u=x_{u^*}$ for any $u\in U$ with $d_U(u)\ge  2$, and $x_w=x_{u^*}$ for any  $w\in W$ with $d_U(w)\ge  1$ {in view of (\ref{eqa})}. Let $w_1,w_2\in W$ and {$w_1w_2\in E(G^*)$}. By Lemma \ref{le-zhai}, 
we have  $d_U(w_1)\ge 1$ and $d_U(w_2)\ge 1$. Since {$G^*$} is $V_5$-free, we get $|N_{G^*}(w)\cap V(H)|\le 1$ for any $w\in W$ and for any $H\in \mathcal{H}$.  Then we get 
 \[ \lambda x_{w_1} 
 = \sum_{u\in N_U(w_1)} x_u + x_{w_2} 
 \le \sum\limits_{u\in U}x_u-2x_u^*+x_{w_2}<\lambda x_{u^*} ,\]
 which leads to $x_{w_1} < x_{u^*}$,   a contradiction. 
 Therefore, we have proved Claim \ref{eq-W0}. 
 \QED

\begin{claim}\label{claa}
There exists a component $H\in \mathcal{H}$ such that $H$ is isomorphic to a star.
\end{claim}

\noindent \textbf{Proof}. Suppose on the contrary that every non-trivial connected component of $G^*[U]$ is isomorphic to $K_3$. Let $t$ be the number of copies of $K_3$ in  $\mathcal{H}$. Then 
{\[ G^*[U]=tK_3\cup (|U|-3t)K_1. \]}
Note that $t\ge 1$ by Claim \ref{claim-eU}. Firstly, we claim that $W\neq \varnothing$. Otherwise, we have {$G^*= R_{t,|U|-3t}$}. Then by Lemma \ref{le-zhai-1}, we get $\lambda< \frac{1+\sqrt{4m-3}}{2}$, a contradiction.

For each non-trivial component $H\in \mathcal{H}$, 
we denote  
\[ W_H=\bigcup_{u\in V(H)}N_W(u). \] 
Let $W_{U_0}$ be the set of vertices of $W$ which are adjacent to a vertex of $U_0$. Namely, 
 \[ W_{U_0}=\bigcup_{u\in U_0}N_W(u). \] 
In view of (\ref{eqc}) and Claim \ref{eq-W0}, 
we obtain
 \begin{align}
   (e(U)+e(U,W)-1)x_{u^*}&\le  \sum\limits_{u\in U}(d_U(u)-1)x_u+\sum\limits_{w\in W}d_U(w)x_w\nonumber\\
                         &=\sum\limits_{H\in\mathcal{H}}\sum\limits_{u\in V(H)}x_u-\sum\limits_{u\in U_0}x_u + 
                         \sum\limits_{w\in W}d_U(w)x_w\nonumber\\
                         &\le \sum\limits_{H\in\mathcal{H}}\left( 
                         \sum\limits_{u\in V(H)}x_u+\sum\limits_{w\in W_H}d_H(w)x_w \right)+\sum\limits_{w\in W_{U_0}}d_{U_0}(w)x_w.\label{eqm}
   \end{align}
Observe that 
\begin{align} 
\sum\limits_{w\in W_{U_0}}d_{U_0}(w)x_w\le e(U_0,W)x_{u^*}\label{eqn}
\end{align} 
 and
 \begin{align}
 e(U)+e(U,W)=\sum\limits_{H\in\mathcal{H}}(e(H)+e(V(H),W))+e(U_0,W).\label{eqo}
 \end{align}
{Combining (\ref{eqm}), (\ref{eqn}) and (\ref{eqo}), we have}
\begin{align}\label{eqd}
 \sum\limits_{H\in\mathcal{H}} 
  \left(e(H)+e(V(H),W) \right) x_{u^*} - x_{u^*} 
  \le 
  \sum\limits_{H\in\mathcal{H}} 
  \left(\sum\limits_{u\in V(H)}x_u+\sum\limits_{w\in W_H}d_H(w)x_w \right). 
\end{align}
Since each component $H\in\mathcal{H}$  is isomorphic to {a} triangle,  we have 
\[ \sum\limits_{u\in V(H)}x_u+\sum\limits_{w\in W_H}d_H(w)x_w\le (e(H)+e(V(H),W))x_{u^*}. \] 
To show a contradiction with (\ref{eqd}), 
 it suffices to prove that there exists a component
  $H\in \mathcal{H}$ satisfying the following strict inequality: 
\[  \sum\limits_{u\in V(H)}x_u+\sum\limits_{w\in W_H}d_H(w)x_w< (e(H)+e(V(H),W)-1)x_{u^*}. \] 

To begin with, 
we fix any one triangle $H\in \mathcal{H}$ with $V(H)=\{u_1,u_2,u_3\}$. 
Next, we start with a discussion of two simple
situations. 

{\bf Case 1.}
 If $H$ contains at least two vertices adjacent to some vertices of $W$, then we can assume that $N_W(u_1)\neq \varnothing $ and $N_W(u_2) \neq \varnothing $. Since $G^*$ is $V_5$-free, $N_W(u_1)\cap N_W(u_2)=\varnothing $. Suppose that $x_{u_1}\ge  x_{u_2}$.  In this case, we can move the edges from $u_2$ to $u_1$, i.e., 
we construct a graph $G'$  such that $V(G')=V(G^*)$ and
\[ E(G')=E(G^*)-\{u_2w: w\in N_W(u_2)\} 
+\{u_1w : w\in N_W(u_2)\}. \] 
 Obviously, the graph $G'$ is $V_5$-free. 
 By  Lemma \ref{le-niki}, we get $\lambda (G')> \lambda $, a contradiction.

{\bf Case 2.} 
If $H$ contains no vertex adjacent to vertices of $W$, then $d_H(w)=0$ for any $w\in W$. 
 By symmetry, we know that $x_{u_1}=x_{u_2}=x_{u_3}$. 
 By $\lambda x_{u_1}=x_{u^*}+2x_{u_1}$, we get 
 $x_{u_1} \le \frac{1}{\lambda -2} x_{u^*}$.
 Since $W\neq \varnothing$ and each vertex of 
 $W$ has at least two neighbors in $U$ by Lemma \ref{le-zhai} and Claim 2, we must have $m\geq 10$. {Recall that $\lambda\ge \frac{1+\sqrt{4m-3}}{2}>\frac{7}{2}$}, we have $\frac{1}{\lambda -2} < \frac{2}{3}$. 
 It follows that $x_{u_1}< \frac{2}{3}x_{u^*}$ and
\begin{align*}
 \sum\limits_{u\in V(H)}x_u+\sum\limits_{w\in W_H}d_H(w)x_w &=3x_{u_1}+\sum\limits_{w\in W_H}d_H(w)x_w  \\
 &<2 x_{u^*}+\sum\limits_{w\in W_H}d_H(w)x_w\\
 &\leq (e(H)+e(V(H),W)-1)x_{u^*}.
\end{align*} 
So $H$ is the desired component contradicting with (\ref{eqd}). 

{\bf Case 3.} 
From the above discussion, we can assume that 
every component $H$ contains exactly one vertex adjacent to some vertices of $W$.  
 Without loss of generality, we may assume that  $N_W(u_1)\neq\varnothing $ and $N_W(u_2)=N_W(u_3)=\varnothing $. 
 Then $x_{u_2}=x_{u_3}$. 
 First of all, we assume that $m\ge 22$.  
Using $\lambda x_{u_2}=x_{u^*}+x_{u_1}+x_{u_2}$,  
we get $x_{u_2}=x_{u_3} \le \frac{2}{\lambda -1} x_{u^*}  
  < \frac{1}{2} x_{u^*}$ as $m\ge  22$. Thus, we have 
\begin{align*}
 \sum\limits_{u\in V(H)}x_u+\sum\limits_{w\in W_H}d_H(w)x_w 
<2 x_{u^*}+\sum\limits_{w\in W_H}d_H(w)x_w 
\le (e(H)+e(V(H),W)-1)x_{u^*}.
\end{align*}
Thus, we conclude that $H$ is a desired component 
contradicting with (\ref{eqd}). 
Next, we suppose that $8\le m\le 21$.  
In this case, we know that $G^*[U]$ contains at most three triangles and some isolated vertices.

{\bf Subcase 3.1.} 
If $G^*[U]$ contains at least two triangles, 
then we have $m\geq 14$ using $W\neq \varnothing$ and 
Lemma \ref{le-zhai}. Let $H_1$ and $H_2$  be two non-trivial components of $G^*[U]$, with  $V(H_1)=\{u_1,u_2,u_3\}$ and $V(H_2)=
\{v_1,v_2,v_3\}$, where $d_{W}(u_1)\ge 1$ and $d_{W}(v_1)\ge 1$. 
By symmetry, we have $x_{u_2}=x_{u_3}=x_{v_2}=x_{v_3}$. 
Using $\lambda x_{u_2}=x_{u^*}+x_{u_1}+x_{u_2}$, we get $x_{u_2}=x_{u_3}=x_{v_2}=x_{v_3}\le \frac{2}{\lambda -1} x_{u^*}  
  < 0.64x_{u^*}$ as $m\ge 14$. One can verify that 
\begin{align*}
\sum\limits_{i=1}^2 
  \left(\sum\limits_{u\in V(H_i)}x_u+\sum\limits_{w\in W_{H_i}}d_{H_i}(w)x_w \right)<
 \left(\sum\limits_{i=1}^2 
  e(H_i)+e(V(H_i),W) \right) x_{u^*} - x_{u^*}. 
\end{align*}
Moreover, we get 
\begin{align*}
\sum\limits_{H\in \mathcal{H}} 
  \left(\sum\limits_{u\in V(H)}x_u+\sum\limits_{w\in W_H}d_H(w)x_w \right)<
 \left(\sum\limits_{H\in \mathcal{H}} 
  e(H)+e(V(H),W) \right) x_{u^*} - x_{u^*}. 
\end{align*}
This is a contradiction to (\ref{eqd}).

\begin{figure}
\centering
\begin{subfigure}[t]{0.2\textwidth}
\centering
 \begin{tikzpicture}
\draw[black](-0.6,0) circle(0.05);
\draw[black](0.6,0) circle(0.05);
\draw (-0.55,0)--(0.55,0);
\draw[black](0,-0.6) circle(0.05);
\draw[black](0,0.6) circle(0.05);
\draw (-0.55,0)--(0,-0.55);
\draw (-0.55,0)--(0,0.55);
\draw (0.55,0)--(0,-0.55);
\draw (0.55,0)--(0,0.55);
\draw (0,-0.55)--(0,0.55);
\draw[black](1.2,0) circle(0.05);
\draw[black](1.6,0) circle(0.05);
\draw (0,0.55)--(1.15,0);
\draw (0,0.55)--(1.55,0);
\draw[black](0.6,-1) circle(0.05);
\draw (0.6,-0.95)--(0,-0.55);
\draw (0.6,-0.95)--(1.15,0);
\draw (0.6,-0.95)--(1.55,0);
\end{tikzpicture}
\subcaption*{$G_1$}
\end{subfigure}
\begin{subfigure}[t]{0.2\textwidth}
\centering
 \begin{tikzpicture}
\draw[black](-0.6,0) circle(0.05);
\draw[black](0.6,0) circle(0.05);
\draw (-0.55,0)--(0.55,0);
\draw[black](0,-0.6) circle(0.05);
\draw[black](0,0.6) circle(0.05);
\draw (-0.55,0)--(0,-0.55);
\draw (-0.55,0)--(0,0.55);
\draw (0.55,0)--(0,-0.55);
\draw (0.55,0)--(0,0.55);
\draw (0,-0.55)--(0,0.55);
\draw[black](1.2,0) circle(0.05);
\draw[black](1.6,0) circle(0.05);
\draw[black](2,0) circle(0.05);
\draw (0,0.55)--(1.15,0);
\draw (0,0.55)--(1.55,0);
\draw (0,0.55)--(1.95,0);
\draw[black](0.6,-1) circle(0.05);
\draw (0.6,-0.95)--(0,-0.55);
\draw (0.6,-0.95)--(1.15,0);
\draw (0.6,-0.95)--(1.55,0);
\end{tikzpicture}
\subcaption*{$G_2$}
\vspace{-5mm}
\end{subfigure}
\caption{The graphs $G_1$ and $G_2$. }
\label{fig-2-graphs}
\end{figure}
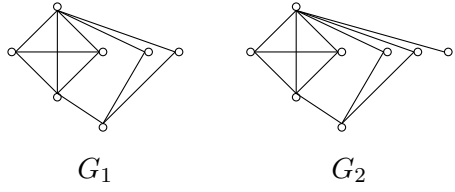

{\bf Subcase 3.2.} 
If $G^*[U]$ contains exactly one triangle, 
then by $W\neq \varnothing$ and Lemma \ref{le-zhai}, 
we know that $|U_0|\ge 1$ and $m\ge 9$. 
Let $H$ be the unique triangle of $G^*[U]$. Let $V(H)=\{u_1,u_2,u_3\}$ and $d_{W}(u_1)\ge 1$. Hence, we have $x_{u_2}=x_{u_3}$. Using $\lambda x_{u_2}=x_{u^*}+x_{u_1}+x_{u_2}$, we get $x_{u_2}=x_{u_3}\le \frac{2}{\lambda -1} x_{u^*}  
  < 0.85 x_{u^*}$ as $m\ge 9$.
If there exists a vertex $w\in W$ with $d_{G^*}(w)=2$, then $x_{w}\le \frac{2}{\lambda} x_{u^*}  
  < 0.6 x_{u^*}$ as $m\ge 9$. 
  However, in this case, it is easy to check that
 \begin{align*}
     \sum\limits_{u\in U}(d_U(u)-1)x_u+\sum\limits_{w\in W}d_U(w)x_w 
     &< 1.7x_{u^*} + x_{u_1} + 1.2x_{u^*} + 
     (e(U,W) -2)x_{u^*} \\
     &<(e(U)+e(U,W)-1)x_{u^*}. 
   \end{align*}
 a contradiction to (\ref{eqm}).
 Hence, there exists no vertex of $W$ with degree two. Combining with Lemma \ref{le-zhai}, for any $w\in W$, we have $d_{G^*}(w)\ge 3$. It follows that $|U_0|\ge 2$ and $m\ge 11$ from $d_W(u_2)=d_W(u_3)=0$. 
 If $m=11$ or $12$, then {$G^*= G_1$ or $G_2$}, respectively; {see Figure \ref{fig-2-graphs}}. By simple calculation, we obtain
$\lambda(G_1)\approx 3.408 <\frac{1+\sqrt{4m-3}}{2} \approx 3.701$ and $\lambda(G_2) \approx 3.487 <\frac{1+\sqrt{4m-3}}{2} \approx 3.854$ , a contradiction.  
If $m\ge 13$, then we get $x_{u_2}=x_{u_3}\le \frac{2}{\lambda -1} x_{u^*}  
  < 0.67 x_{u^*}$. Since $ V_5 \subseteq K_5$ 
  and $G^*$ is $K_5$-free, we get 
  from (\ref{eq1}) that $\lambda \le \sqrt{(1-1/4)2m}$. Thus, 
  for any $u\in U_0$, we obtain $x_u\ge \frac{1}{\lambda}x_{u^*} \ge \frac{1}{\sqrt{(1-1/4)2m}}x_{u^*} \ge 0.178 x_{u^*}$ as $m\le 21$. It is easy to verify that 
  \begin{align*}
    \sum\limits_{u\in U}(d_U(u)-1)x_u+\sum\limits_{w\in W}d_U(w)x_w 
    < (e(U)+e(U,W)-1)x_{u^*}, 
   \end{align*}
which is a contradiction to (\ref{eqm}). 

Therefore, we have deduced a contradiction in 
all cases. 
We conclude that there exists at least one non-trivial connected component isomorphic to a star in $G^*[U]$.\QED

\medskip 

From Claim \ref{claa}, 
$G^*[U]$ contains at least one non-trivial connected component isomorphic to a star. From the definition in (\ref{eq-eta}), we have 
 $\eta_1(K_3)\le 0$ and $\eta_1(K_{1,s})\le -1$  for any $s\ge  1$. 
 It follows from 
 (\ref{eqb}) that there is exactly one non-trivial connected component isomorphic to a star in $G^*[U]$. 
 Let {$H' = K_{1,s}$} be such a star  for some $s\ge  1$, 
 and let $t$  be number of triangles in $G^*[U]$. 
 Note that the equality in (\ref{eqb}) holds, which implies 
the following: 
\begin{itemize}

\item[(a)]  
 $e(W)=0$, 
 $U_0=\varnothing $ and $G^*[U] = t K_3 \cup K_{1,s}$  
 for $s\ge 1$ and $t\ge 0$;

 \item[(b)] 
 $x_u=x_{u^*}$ for any $u\in U$ with $d_U(u)\ge  2$; 
 
 \item[(c)] 
 $x_w=x_{u^*}$ for any  $w\in W$ with $d_U(w)\ge  1$.  
 \end{itemize}

\begin{claim}\label{cl-we}
$W=\varnothing $.
\end{claim}

\noindent \textbf{Proof.} 
Suppose on the contrary that $W\neq \varnothing $.  
We fix a vertex $w\in W$. 
Let $V(H')=\{v_0,v_1,\ldots,v_s\}$, where $v_0$ is the central vertex of $H'$. 
If $s\ge  2$, then  $w$ can not be adjacent to all vertices of $\{v_0,v_1,v_2\}$. Since $e(W)=0$, we have $\lambda x_w =
\sum_{u\in N_U(w)} x_u <  \sum_{u\in U} x_u = \lambda x_{u^*}$, which leads to $x_w< x_{u^*}$. This contradicts with (c). 
If $s=1$ and $t=  0$, then it follows from {(a)} that 
$\lambda x_{u^*} =\sum_{u\in U} x_u
\le 2x_{u^*}$, which implies that $\lambda \le 2<\frac{1+\sqrt{4m-3}}{2}$ when  $m\ge  8$, a contradiction. 
If $s=1$ and $t \ge 1$, then 
the vertex $w$ can not be adjacent to all vertices of a  
copy of $K_3$ in $G^*[U]$ since $G^*$ is $V_5$-free. 
Similarly, {we have $\lambda x_w =
\sum_{u\in N_U(w)} x_u <  \sum_{u\in U} x_u = \lambda x_{u^*}$, which leads to $x_w< x_{u^*}$.
It is a contradiction  to (c).}
with (c).\QED

\medskip 
To sum up, we know that $G^*= K_1\vee (K_{1,s}\cup tK_3)$. If $t \neq 0$, then $x_u= x_{u^*}$ does not hold for any vertex $u$ in the component of $G^*[U]$ which is isomorphic to $K_3$, contradicting with (b).
 Thus, we get $t=0$ and so 
 $G^*=K_1\vee K_{1,s}$. Observe that 
 $2s+1=m$. Then $G^*=K_2\vee \frac{m-1}{2}K_1$. 
 We complete the proof of Theorem \ref{Thm-V5}.

\section{Proof of Theorem \ref{Thm-F3}}

In this section, we will consider the maximum spectral radius of $F_3$-free graphs with $m$ edges. 
{In this case, $G^*$ is a spectral extremal graph for the  graph $F_3$, and $u^*$ is an extremal vertex of $G^*$.}  Let $\lambda (G^*)=\lambda$ and $\bf{x}$ be the Perron vector of $G^*$ with 
coordinate $x_u$ corresponding to the vertex $u\in V(G)$. 
For simplicity, let $U=N(u^*)$ and 
$W=V(G^*)\setminus N[u^*]$. 
Since $G^*$ is $F_3$-free, we know that 
$G^*[U]$ does not contain a matching consisting of 
three edges.

\begin{claim}
$G^*$ is connected.
\end{claim}
\noindent \textbf{Proof.} 
Otherwise, let $H_1$ and $H_2$ be two non-trivial components of $G^*$ with $\lambda (H_1)=\lambda (G^*)$. Choosing a vertex $v\in V(H_1)$ and an edge $w_1w_2\in E(H_2)$, we construct a new graph {$G'=G^*-w_1w_2+vw_1$}. Clearly $G'$ is still $F_3$-free and $\lambda (G^*) <\lambda (G')$, a contradiction.\QED

Since $G^*$ is connected, $\mathbf{x}$ is  a positive real vector. Note that $K_3\vee \frac{m-3}{3}K_1$ is 
$F_3$-free. Then 
$ \lambda \ge  \lambda (K_3\vee \frac{m-3}{3}K_1)=1+\sqrt{m-2}$. 
Hence, we get $\lambda ^{2}-2\lambda \ge  m-3$. 
Recall that 
\begin{align*}
\lambda^{2}x_{u^*}=|U|x_{u^*}+\sum\limits_{u\in U}d_U(u)x_u+\sum\limits_{w\in W}d_U(w)x_w. 
\end{align*}
{Combining with $\lambda =\sum_{u\in U}x_u$,} we obtain
\begin{align*}
    (\lambda ^{2}-2\lambda )x_{u^*}=|U|x_{u^*}+\sum\limits_{u\in U}d_U(u)x_u+\sum\limits_{w\in W}d_U(w)x_w-\sum\limits_{u\in U}2x_u.
\end{align*}
Notice that $\lambda ^{2}-2\lambda \ge  m-3=|U|+e(U)+e(W)+e(U,W)-3$. Then
\begin{align*}
   (e(U)+e(W)+e(U,W)-3)x_{u^*}\le \sum\limits_{u\in U}{(d_U(u)-2)}x_u+\sum\limits_{w\in W}d_U(w)x_w.
\end{align*}
By simplifying, we have
\begin{align}
   e(W)&\le \sum\limits_{u\in U}{(d_U(u)-2)}\frac{x_u}{x_{u^*}}-e(U)+\sum\limits_{w\in W}d_U(w)\frac{x_w}{x_{u^*}}-e(U,W)+3\notag\\
   &\le  \sum\limits_{u\in U}(d_U(u)-2)\frac{x_u}{x_{u^*}}-e(U)+3,\label{eqi}
\end{align}
 and the equality occurs if and only if $\lambda ^{2}-2\lambda = m-3$ and $x_w=x_{u^*}$ holds for each vertex $w\in W$ with $d_U(w)\ge 1$.

 \begin{claim}\label{cl-ubo}
$\lambda \le \frac{1}{2}(e(U)+3)$.
\end{claim}
 \noindent \textbf{Proof.} 
 Since $\lambda x_{u^*} = \sum_{u\in U} x_u$  
 and $e(W)\ge  0$, it follows from (\ref{eqi}) that 
 \begin{align*}
 2 \lambda &= 2 \sum\limits_{u\in U}\frac{x_u}{x_{u^*}}\le \sum\limits_{u\in U}d_U(u)\frac{x_u}{x_{u^*}}-e(U)+3 \le e(U)+3. 
\end{align*}
The proof is completed.\QED

\begin{claim}\label{cl-com1}
$G^*[U]$ contains one or two non-trivial components.
\end{claim}

\noindent \textbf{Proof.}
Observe that  $\lambda = 
1+\sqrt{m-2}>\sqrt{m}$. We have 
\begin{align*}
mx_{u^*}&<\lambda^{2}x_{u^*} 
=|U|x_{u^*}+\sum\limits_{u\in U} d_U(u)x_u+{\sum\limits_{w\in W} d_U(w)x_w}\\
&\le |U|x_{u^*}+2e(U)x_{u^*}+e(U,W)x_{u^*}.
\end{align*}
Notice that $m=|U|+e(U)+e(W)+e(U,W)$. 
Thus, we get $e(W)<e(U)$. 
If $e(U)=0$, then $e(W)<0$, a contradiction.
Hence, $G^*[U]$ contains at least one non-trivial {component}. Since $G^*$ is $F_3$-free, 
there are not three non-trivial components in $G^*[U]$. Therefore, $G^*[U]$ contains one or two non-trivial components, as desired. 
\QED

\medskip 
  By Claim \ref{cl-com1}, $G^*[U]$ contains one or two non-trivial components. Let $\mathcal{H}$ be the set of the non-trivial components of $G^*[U]$. 
  For each $H\in \mathcal{H}$, we denote 
  \[ \eta_2(H):=\sum\limits_{u\in V(H)}(d_H(u)-2)\frac{x_u}{x_{u^*}}-e(H). \]
  Let $U_0=\{u\in U|d_U(u)=0\}$. 
  Then (\ref{eqi}) can be rewritten as
\begin{align}\label{eqg}
e(W)\le \sum\limits_{H\in \mathcal{H}}\eta_2(H)-2\sum\limits_{u\in U_0}\frac{x_u}{x_{u^*}}+3.
\end{align}
 Furthermore, since $e(W)\ge 0$, we can get
\begin{align}\label{eqf}
\sum\limits_{H\in \mathcal{H}}\eta_2(H)\ge  e(W)+ 2\sum\limits_{u\in U_0}\frac{x_u}{x_{u^*}}-3\ge -3,
\end{align}
 and the last equality holds if and only if $e(W)=0$ and $U_0=\varnothing $.

\begin{claim}\label{cl-comp}
$G^*[U]$ contains exactly one non-trivial connected component.
\end{claim}
\noindent \textbf{Proof.} 
According to Claim \ref{cl-com1}, we assume on the contrary that 
$G^*[U]$ contains two non-trivial connected components. As $G^*[U]$ is $3K_2$-free,  each non-trivial component $H\in \mathcal{H}$ 
is $2K_2$-free. 
Thus, $H$ is isomorphic to either triangle or star. Note that $\eta_2(K_3)=-3$ and $\eta_2(K_{1,t})<-2$ for any $t\ge  1$. It follows that  $\sum\limits_{H\in \mathcal{H}}\eta_2(H)<-4$, contradicting with (\ref{eqf}). \QED

By Claim \ref{cl-comp}, there is exactly one non-trivial connected component. We denote it by $H^*$. Then (\ref{eqg}) can be transformed to
\begin{align}\label{eqh}
e(W)\le \eta_2(H^*)-2\sum\limits_{u\in U_0}\frac{x_u}{x_{u^*}}+3.
\end{align}
 And we have to emphasize that from (\ref{eqf}), 
 we get 
 \begin{align}\label{eql}
     \eta_2(H^*)\ge  -3 
 \end{align}
 and the equality holds if and only if $\lambda ^{2}-2\lambda = m-3$, $x_w=x_{u^*}$ holds for any $w\in W$ with $d_U(w)\ge 1$,  $e(W)=0$  and $U_0=\varnothing$.

\begin{claim}\label{cl-mind}
$\delta(H^*)\ge  2$ for $m\ge  28$.
\end{claim}
\noindent \textbf{Proof.} Otherwise, 
suppose on the contrary that 
$\delta(H^*)=1$. Note that $G^*$ is $F_3$-free, $H^*$ is $3K_2$-free. From Lemma \ref{le-Liu}, we know that $H^*$ contains a dominating complete bipartite subgraph. Choose a maximal dominating complete bipartite subgraph and denote the two color classes of this bipartite subgraph by $S$ and $T$, where $|S|=s$ and $|T|=t$. Without loss of generality, assume that $1\le s\le t$. 
Since $H^*$ is $3K_2$-free, 
we get $s\le 2$ immediately. In what follows, we denote  $X=S\cup T$ and $Y=V(H^*)\setminus X$.

{\bf Case 1}. $s=2$.
  
  Since $\delta(H^*)=1$ and $d_{H^*}(u)\ge  2$ for any $u\in X$,  we know that $Y$ must contain a vertex with degree one in $H^*$. Notice that  $H^*$ is $3K_2$-free, so there is no edge in {$H^*[Y]$}. 
  
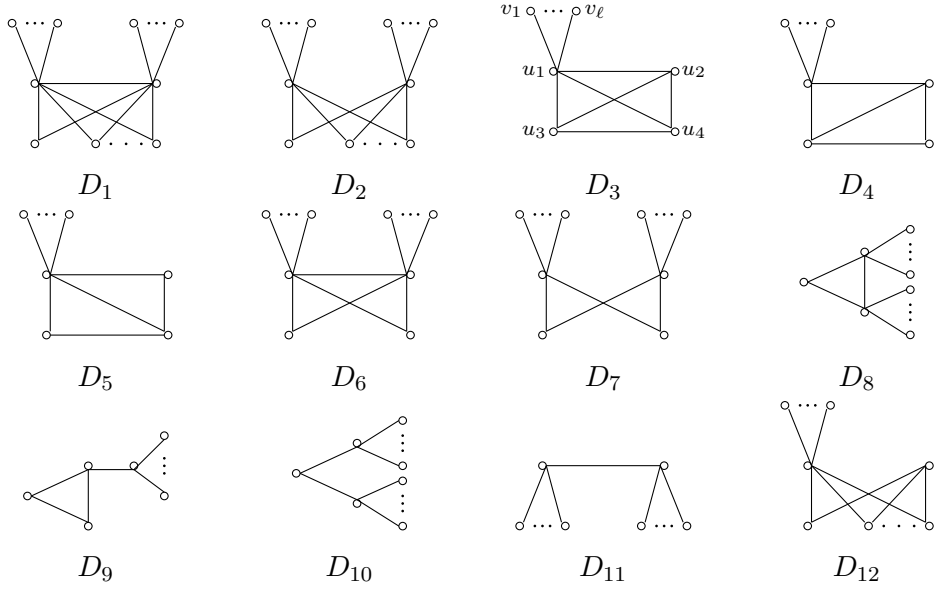
\begin{figure}
\centering
\begin{subfigure}[t]{0.2\textwidth}
\centering
 \begin{tikzpicture}
\draw[black](-0.8,0) circle(0.05);
\draw[black](0.8,0) circle(0.05);
\draw (-0.75,0)--(0.75,0);
\draw[black](0,-0.8) circle(0.05);
\draw (-0.75,0)--(-0.05,-0.75);
\draw (0.75,0)--(0.05,-0.75);
\draw[black](-0.8,-0.8) circle(0.05);
\draw (-0.75,0)--(-0.75,-0.75);
\draw (0.75,0)--(-0.75,-0.75);
\draw[black](0.8,-0.8) circle(0.05);
\draw (-0.75,0)--(0.75,-0.75);
\draw (0.75,0)--(0.75,-0.75);
\draw[black](0.5,0.8) circle(0.05);
\draw[black](1.1,0.8) circle(0.05);
\fill(0.8, 0.8) circle(0.5pt); 
\fill(0.7, 0.8) circle(0.5pt); 
\fill(0.9, 0.8) circle(0.5pt); 
\draw[black](-0.5,0.8) circle(0.05);
\draw[black](-1.1,0.8) circle(0.05);
\fill(-0.8, 0.8) circle(0.5pt); 
\fill(-0.7, 0.8) circle(0.5pt); 
\fill(-0.9, 0.8) circle(0.5pt); 

\draw (-0.75,0)--(-0.55,0.75);
\draw (-0.75,0)--(-1.05,0.75);
\draw (0.75,0)--(0.55,0.75);
\draw (0.75,0)--(1.05,0.75);
\fill(0.2, -0.8) circle(0.5pt); 
\fill(0.4, -0.8) circle(0.5pt); 
\fill(0.6, -0.8) circle(0.5pt); 

\end{tikzpicture}
\subcaption*{$D_1$}
\end{subfigure}
\begin{subfigure}[t]{0.2\textwidth}
\centering
 \begin{tikzpicture}
\draw[black](-0.8,0) circle(0.05);
\draw[black](0.8,0) circle(0.05);
\draw[black](0,-0.8) circle(0.05);
\draw (-0.75,0)--(-0.05,-0.75);
\draw (0.75,0)--(0.05,-0.75);
\draw[black](-0.8,-0.8) circle(0.05);
\draw (-0.75,0)--(-0.75,-0.75);
\draw (0.75,0)--(-0.75,-0.75);
\draw[black](0.8,-0.8) circle(0.05);
\draw (-0.75,0)--(0.75,-0.75);
\draw (0.75,0)--(0.75,-0.75);
\draw[black](0.5,0.8) circle(0.05);
\draw[black](1.1,0.8) circle(0.05);
\draw[black](-0.5,0.8) circle(0.05);
\draw[black](-1.1,0.8) circle(0.05);
\draw (-0.75,0)--(-0.55,0.75);
\draw (-0.75,0)--(-1.05,0.75);
\draw (0.75,0)--(0.55,0.75);
\draw (0.75,0)--(1.05,0.75);
\fill(0.8, 0.8) circle(0.5pt); 
\fill(0.7, 0.8) circle(0.5pt); 
\fill(0.9, 0.8) circle(0.5pt); 
\fill(-0.8, 0.8) circle(0.5pt); 
\fill(-0.7, 0.8) circle(0.5pt); 
\fill(-0.9, 0.8) circle(0.5pt); 
\fill(0.2, -0.8) circle(0.5pt); 
\fill(0.4, -0.8) circle(0.5pt); 
\fill(0.6, -0.8) circle(0.5pt); 
\end{tikzpicture}
\subcaption*{$D_2$}
\vspace{-5mm}
\end{subfigure}
\begin{subfigure}[t]{0.2\textwidth}
\centering
 \begin{tikzpicture}
\draw[black](-0.8,0) circle(0.05);
\draw(-1.05,0) node{{\scriptsize $u_1$}};

\draw[black](0.8,0) circle(0.05);
\draw(1.05,0) node{{\scriptsize $u_2$}};
\draw (-0.75,0)--(0.75,0);
\draw[black](-0.8,-0.8) circle(0.05);
\draw(-1.05,-0.8) node{{\scriptsize $u_3$}};
\draw (-0.75,0)--(-0.75,-0.75);
\draw (0.75,0)--(-0.75,-0.75);
\draw (-0.75,-0.8)--(0.75,-0.8);
\draw[black](0.8,-0.8) circle(0.05);
\draw(1.05,-0.8) node{{\scriptsize $u_4$}};
\draw (-0.75,0)--(0.75,-0.75);
\draw (0.75,0)--(0.75,-0.75);
\draw[black](-0.5,0.8) circle(0.05);
\draw(-0.25,0.8) node{{\scriptsize $v_{\ell}$}};
\draw[black](-1.1,0.8) circle(0.05);
\draw(-1.35,0.8) node{{\scriptsize $v_1$}};
\draw (-0.75,0)--(-0.55,0.75);
\draw (-0.75,0)--(-1.05,0.75);
\fill(-0.8, 0.8) circle(0.5pt); 
\fill(-0.7, 0.8) circle(0.5pt); 
\fill(-0.9, 0.8) circle(0.5pt); 

\end{tikzpicture}
\subcaption*{$D_3$}
\end{subfigure}
\begin{subfigure}[t]{0.2\textwidth}
\centering
 \begin{tikzpicture}
\draw[black](-0.8,0) circle(0.05);
\draw[black](0.8,0) circle(0.05);
\draw (-0.75,0)--(0.75,0);
\draw[black](-0.8,-0.8) circle(0.05);
\draw (-0.75,0)--(-0.75,-0.75);
\draw (0.75,0)--(-0.75,-0.75);
\draw (-0.75,-0.8)--(0.75,-0.8);
\draw[black](0.8,-0.8) circle(0.05);
\draw (0.75,0)--(0.75,-0.75);
\draw[black](-0.5,0.8) circle(0.05);
\draw[black](-1.1,0.8) circle(0.05);
\draw (-0.75,0)--(-0.55,0.75);
\draw (-0.75,0)--(-1.05,0.75);
\fill(-0.8, 0.8) circle(0.5pt); 
\fill(-0.7, 0.8) circle(0.5pt); 
\fill(-0.9, 0.8) circle(0.5pt); 
\end{tikzpicture}
\subcaption*{$D_4$}
\end{subfigure}
\begin{subfigure}[t]{0.2\textwidth}
\centering
 \begin{tikzpicture}
\draw[black](-0.8,0) circle(0.05);
\draw[black](0.8,0) circle(0.05);
\draw (-0.75,0)--(0.75,0);
\draw[black](-0.8,-0.8) circle(0.05);
\draw (-0.75,0)--(-0.75,-0.75);
\draw (-0.75,0)--(0.75,-0.75);
\draw (-0.75,-0.8)--(0.75,-0.8);
\draw[black](0.8,-0.8) circle(0.05);
\draw (0.75,0)--(0.75,-0.75);
\draw[black](-0.5,0.8) circle(0.05);
\draw[black](-1.1,0.8) circle(0.05);
\draw (-0.75,0)--(-0.55,0.75);
\draw (-0.75,0)--(-1.05,0.75);
\fill(-0.8, 0.8) circle(0.5pt); 
\fill(-0.7, 0.8) circle(0.5pt); 
\fill(-0.9, 0.8) circle(0.5pt); 
\end{tikzpicture}
\subcaption*{$D_5$}
\end{subfigure}
\begin{subfigure}[t]{0.2\textwidth}
\centering
 \begin{tikzpicture}
\draw[black](-0.8,0) circle(0.05);
\draw[black](0.8,0) circle(0.05);
\draw (-0.75,0)--(0.75,0);

\draw[black](-0.8,-0.8) circle(0.05);
\draw (-0.75,0)--(-0.75,-0.75);
\draw (0.75,0)--(-0.75,-0.75);
\draw[black](0.8,-0.8) circle(0.05);
\draw (-0.75,0)--(0.75,-0.75);
\draw (0.75,0)--(0.75,-0.75);
\draw[black](0.5,0.8) circle(0.05);
\draw[black](1.1,0.8) circle(0.05);
\draw[black](-0.5,0.8) circle(0.05);
\draw[black](-1.1,0.8) circle(0.05);
\draw (-0.75,0)--(-0.55,0.75);
\draw (-0.75,0)--(-1.05,0.75);
\draw (0.75,0)--(0.55,0.75);
\draw (0.75,0)--(1.05,0.75);
\fill(0.8, 0.8) circle(0.5pt); 
\fill(0.7, 0.8) circle(0.5pt); 
\fill(0.9, 0.8) circle(0.5pt); 
\fill(-0.8, 0.8) circle(0.5pt); 
\fill(-0.7, 0.8) circle(0.5pt); 
\fill(-0.9, 0.8) circle(0.5pt); 
\end{tikzpicture}
\subcaption*{$D_6$}
\end{subfigure}
\begin{subfigure}[t]{0.2\textwidth}
\centering
 \begin{tikzpicture}
\draw[black](-0.8,0) circle(0.05);
\draw[black](0.8,0) circle(0.05);
\draw[black](-0.8,-0.8) circle(0.05);
\draw (-0.75,0)--(-0.75,-0.75);
\draw (0.75,0)--(-0.75,-0.75);
\draw[black](0.8,-0.8) circle(0.05);
\draw (-0.75,0)--(0.75,-0.75);
\draw (0.75,0)--(0.75,-0.75);
\draw[black](0.5,0.8) circle(0.05);
\draw[black](1.1,0.8) circle(0.05);
\draw[black](-0.5,0.8) circle(0.05);
\draw[black](-1.1,0.8) circle(0.05);
\draw (-0.75,0)--(-0.55,0.75);
\draw (-0.75,0)--(-1.05,0.75);
\draw (0.75,0)--(0.55,0.75);
\draw (0.75,0)--(1.05,0.75);
\fill(0.8, 0.8) circle(0.5pt); 
\fill(0.7, 0.8) circle(0.5pt); 
\fill(0.9, 0.8) circle(0.5pt); 
\fill(-0.8, 0.8) circle(0.5pt); 
\fill(-0.7, 0.8) circle(0.5pt); 
\fill(-0.9, 0.8) circle(0.5pt); 
\end{tikzpicture}
\subcaption*{$D_7$}
\end{subfigure}
\begin{subfigure}[t]{0.2\textwidth}
\centering
 \begin{tikzpicture}
\draw[black](0,0.4) circle(0.05);
\draw[black](0,-0.4) circle(0.05);
\draw (0,-0.35)--(0,0.35);

\draw[black](-0.8,0) circle(0.05);
\draw (-0.75,0)--(0,-0.35);
\draw (-0.75,0)--(0,0.35);
\draw[black](0.6,-0.1) circle(0.05);
\draw[black](0.6,0.1) circle(0.05);
\draw (0.55,-0.1)--(0,-0.35);
\draw (0.55,0.1)--(0,0.35);
\draw[black](0.6,-0.7) circle(0.05);
\draw[black](0.6,0.7) circle(0.05);
\draw (0.55,-0.7)--(0,-0.35);
\draw (0.55,0.7)--(0,0.35);
\fill(0.6, -0.5) circle(0.5pt); 
\fill(0.6, -0.4) circle(0.5pt); 
\fill(0.6, -0.3) circle(0.5pt); 
\fill(0.6, 0.5) circle(0.5pt); 
\fill(0.6, 0.4) circle(0.5pt); 
\fill(0.6, 0.3) circle(0.5pt); 
\end{tikzpicture}
\subcaption*{$D_8$}
\end{subfigure}

\begin{subfigure}[t]{0.2\textwidth}
\centering
 \begin{tikzpicture}
\draw[black](0,0.4) circle(0.05);
\draw[black](0,-0.4) circle(0.05);
\draw (0,-0.35)--(0,0.35);

\draw[black](-0.8,0) circle(0.05);
\draw (-0.75,0)--(0,-0.35);
\draw (-0.75,0)--(0,0.35);
\draw[black](0.6,0.4) circle(0.05);
\draw (0.6,0.35)--(0,0.35);
\draw[black](1,0) circle(0.05);
\draw[black](1,0.8) circle(0.05);
\draw (0.6,0.35)--(1,0.05);
\draw (0.6,0.35)--(1,0.75);
\fill(1, 0.3) circle(0.5pt); 
\fill(1, 0.4) circle(0.5pt); 
\fill(1, 0.5) circle(0.5pt); 
\end{tikzpicture}
\subcaption*{$D_{9}$}
\end{subfigure}
\begin{subfigure}[t]{0.2\textwidth}
\centering
 \begin{tikzpicture}
\draw[black](0,0.4) circle(0.05);
\draw[black](0,-0.4) circle(0.05);

\draw[black](-0.8,0) circle(0.05);
\draw (-0.75,0)--(0,-0.35);
\draw (-0.75,0)--(0,0.35);
\draw[black](0.6,-0.1) circle(0.05);
\draw[black](0.6,0.1) circle(0.05);
\draw (0.55,-0.1)--(0,-0.35);
\draw (0.55,0.1)--(0,0.35);
\draw[black](0.6,-0.7) circle(0.05);
\draw[black](0.6,0.7) circle(0.05);
\draw (0.55,-0.7)--(0,-0.35);
\draw (0.55,0.7)--(0,0.35);
\fill(0.6, -0.5) circle(0.5pt); 
\fill(0.6, -0.4) circle(0.5pt); 
\fill(0.6, -0.3) circle(0.5pt); 
\fill(0.6, 0.5) circle(0.5pt); 
\fill(0.6, 0.4) circle(0.5pt); 
\fill(0.6, 0.3) circle(0.5pt); 
\end{tikzpicture}
\subcaption*{$D_{10}$}
\end{subfigure}
\begin{subfigure}[t]{0.2\textwidth}
\centering
 \begin{tikzpicture}
\draw[black](-0.8,0) circle(0.05);
\draw[black](0.8,0) circle(0.05);
\draw (-0.75,0)--(0.75,0);
\draw[black](0.5,-0.8) circle(0.05);
\draw[black](1.1,-0.8) circle(0.05);
\draw[black](-0.5,-0.8) circle(0.05);
\draw[black](-1.1,-0.8) circle(0.05);
\draw (-0.75,0)--(-0.55,-0.75);
\draw (-0.75,0)--(-1.05,-0.75);
\draw (0.75,0)--(0.55,-0.75);
\draw (0.75,0)--(1.05,-0.75);
\fill(0.9, -0.8) circle(0.5pt); 
\fill(0.8, -0.8) circle(0.5pt); 
\fill(0.7, -0.8) circle(0.5pt); 
\fill(-0.9, -0.8) circle(0.5pt); 
\fill(-0.8, -0.8) circle(0.5pt); 
\fill(-0.7, -0.8) circle(0.5pt); 
\end{tikzpicture}
\subcaption*{$D_{11}$}
\end{subfigure}
\begin{subfigure}[t]{0.2\textwidth}
\centering
 \begin{tikzpicture}
\draw[black](-0.8,0) circle(0.05);
\draw[black](0.8,0) circle(0.05);
\draw[black](0,-0.8) circle(0.05);
\draw (-0.75,0)--(-0.05,-0.75);
\draw (0.75,0)--(0.05,-0.75);
\draw[black](-0.8,-0.8) circle(0.05);
\draw (-0.75,0)--(-0.75,-0.75);
\draw (0.75,0)--(-0.75,-0.75);
\draw[black](0.8,-0.8) circle(0.05);
\draw (-0.75,0)--(0.75,-0.75);
\draw (0.75,0)--(0.75,-0.75);
\draw[black](-0.5,0.8) circle(0.05);
\draw[black](-1.1,0.8) circle(0.05);
\draw (-0.75,0)--(-0.55,0.75);
\draw (-0.75,0)--(-1.05,0.75);
\fill(-0.8, 0.8) circle(0.5pt); 
\fill(-0.7, 0.8) circle(0.5pt); 
\fill(-0.9, 0.8) circle(0.5pt); 
\fill(0.2, -0.8) circle(0.5pt); 
\fill(0.4, -0.8) circle(0.5pt); 
\fill(0.6, -0.8) circle(0.5pt); 

\end{tikzpicture}
\subcaption*{$D_{12}$}
\end{subfigure}


\caption{Forms of $H^*$ when $\delta(H^*)=1$. }
\label{fig-12-graphs}
\end{figure}

  If $t\ge  3$,  then each vertex $y\in Y$ can not be adjacent to a vertex of $T$ as  $H^*$ is $3K_2$-free. {Since $H^*[S,T]$  is a maximal dominating complete subgraph}, $y$ is adjacent to exactly one vertex of  $S$. Hence, $d_{H^*}(y)=1$ holds for each $y\in Y$. Furthermore, 
  $H^*$ is isomorphic to $D_1$ or $D_2$ in Figure \ref{fig-12-graphs}. However, by simple calculation, we have $\eta_2(D_1)<-3$ and $\eta_2(D_2)<-4$, a contradiction to (\ref{eql}). 
  
  If $t=2$,  then because of the maximality of $H^*[S,T]$ and the fact that $Y$ has a vertex of degree one in $H^*$,  we can easily deduce that $d_{H^*}(y)=1$ for any $y\in Y$.  Furthermore, we observe that  $H^*$ is isomorphic to {one of 
  $\{D_3,D_4,D_5,D_6,D_7\}$} in Figure \ref{fig-12-graphs}. By straightforward computations, 
  we know that $\eta_2(D_4)<-3$, $\eta_2(D_5)<-3$, $\eta_2(D_6)<-3$ and $\eta_2(D_7)<-4$. This is a contradiction to (\ref{eql}). In this case, we obtain that 
  $H^*$ is isomorphic to $D_3$.

  In what follows, we deduce a contradiction whenever  
 $H^*$ is isomorphic to $D_3$. 
   Let $u_1$ be the vertex with maximum degree in $H^*$. Let $u_i$ and $v_j$ be the vertex of degree 3 and of degree 1 in $H^*$  for $2\le i\le 4$ and {$1\le j\le \ell$}, respectively; see Figure \ref{fig-12-graphs}. Then in view of  (\ref{eqh}), we know 
  \begin{align*}
    0\le e(W)&\le \eta_2(H^*)-2\sum\limits_{u\in U_0}\frac{x_u}{x_{u^*}}+3\\
    &=\sum\limits_{i=2}^{4}\frac{x_{u_i}}{x_{u^*}}
  +(\ell+1)\frac{x_{u_1}}{x_{u^*}}-\sum\limits_{i=1}^{\ell}\frac{x_{v_i}}{x_{u^*}}-(\ell+6)-2\sum\limits_{u\in U_0}\frac{x_u}{x_{u^*}}+3.  
  \end{align*}
  The above inequality implies that
  \begin{align*}
\sum\limits_{i=1}^{\ell}x_{v_i}+2\sum\limits_{u\in U_0}x_u\le \sum\limits_{i=2}^{4}x_{u_i}
  +(\ell+1)x_{u_1}-(\ell+3)x_{u^*}\le 2x_{u^*}.  
  \end{align*}
  Moreover, we get
  \begin{align*}
   \lambda x_{u^*} 
={\sum\limits_{i=1}^{4}x_{u_i}}+\sum\limits_{i=1}^{\ell}x_{v_i}+\sum\limits_{u\in U_0}x_u\le 6x_{u^*}.    
  \end{align*}
    Hence, we have         $\lambda \le 6$, which is a contradiction since $\lambda \ge 1+\sqrt{m-2}$ and $m\geq 28$.

{\bf Case 2}. $s=1$.  

We may assume that {$S=\{u'\}$}. From the maximality of $H^*[S,T]$, there is no vertex of $Y$ adjacent to {$u'$}. Since $H^*$ is $3K_2$-free and $\delta(H^*)=1$, we obtain that $H^*[T]$ is $P_4$-free.

{\bf Subcase 2.1}.  $H^*[T]$ contains a $P_3$.   Let $v_1v_2v_3 $ be a path of length 2 in $H^*[T]$. Since $H^*$ is $3K_2$-free, $H^*[Y]$ contains no edge.  If $v_1$ is adjacent to  $v_3$, then we can get that $H^*$ has the form of $D_3$. However,  it's impossible as the proof in Case 1 shows. Thus, $v_1$ is not adjacent to $v_3$. {In this case,  if $t=3$, then $H^*$ has the form of   $D_4$ or $D_5$.} It is also impossible as $\eta_2(D_4)<-3$ and $\eta_2(D_5)<-3$. Hence, $t\ge 4$. Then since $H^*$ is $3K_2$-free and $\delta(H^*)=1$, $H[T\setminus\{v_1,v_2,v_3\}]$ contains no edge. Furthermore, we have that $y$ is  only adjacent to $v_2$ for any $y\in Y$. Now  we can get that $H^*$ has the form of {$D_1$. However, 
we have $\eta_2(D_1)<-3$,} a contradiction to (\ref{eql}).

{\bf Subcase 2.2}.  $H^*[T]$ contains no $P_3$, but contains a $P_2$. Since $H^*$ is $3K_2$-free and $\delta(H^*)=1$, $H^*[T]$ contains no $2P_2$.  Therefore, $H^*[T]$ induces exactly one edge. We denote the unique edge by $v_1v_2$.  

If $t\ge  3$, we get that  $H^*[Y]$ contains no edge as $H^*$ is $3K_2$-free.  Besides, we can know that $d_T(y)\le 2$ holds for any $y\in Y$. Otherwise, $H^*$ will contain a $3K_2$.  If there exists a vertex $y'\in Y$ such that $d_T(y')=2$, then $N_T(y')=\{v_1,v_2\}$ as $H^*$ is $3K_2$-free and $\delta(H^*)=1$. Furthermore, $|N_T(y)|=1$ for any {$y\in Y\setminus\{y'\}$}. Hence, $H^*$ has the form of $D_4$. Similarly, we know $\eta_2(D_4)<-3$, and it is contradicted to (\ref{eql}). Therefore,  $|N_T(y)|=1$ for any $y\in Y$ or $Y=\varnothing$. It follows that $H^*$ has the form of $D_8$ or $D_9$. However, {it is  easy} to calculate that $\eta_2(D_8)<-3$ and $\eta_2(D_9)<-4$, a contradiction to (\ref{eql}). So we have $t=2$. 

However, now we will obtain a $K_{1,3}$ which dominates $H^*$. It {contradicts} to the maximality of $H^*[S,T]$.

{\bf Subcase 2.3}. $H^*[T]$ contains  no $P_2$.
If $Y=\varnothing$, then {$H^*= K_{1,t}$}. If $t=1$, then by Claim \ref{cl-ubo}, we have 
$\lambda \le \frac{1}{2}(e(U)+3) =2 <1+\sqrt{m-2}$, a contradiction. 
If $t\ge 2$, then let $v$  be the central vertex of $H^*$, and $v_i$ be the leaf vertex of $H^*$ for each $1\le i\leq t$. 
In view of  (\ref{eqh}), we have
\begin{align*}
 0\le e(W)&\le \eta_2(H^*)-2\sum\limits_{u\in U_0}\frac{x_u}{x_{u^*}}+3\\
    &=(t-2)\frac{x_{v}}{x_{u^*}}-\sum\limits_{i=1}^{t}\frac{x_{v_i}}{x_{u^*}}-t-2\sum\limits_{u\in U_0}\frac{x_u}{x_{u^*}}+3.  
\end{align*}
It follows that 
\begin{align*} {\sum\limits_{i=1}^{t}x_{v_i}+2\sum\limits_{u\in U_0}x_u \le(t-2)x_{v}-tx_{u^*}+3x_{u^*}\le x_{u^*}.}
\end{align*}
Furthermore, we obtain
\begin{align*}
\lambda x_{u^*}= x_v+\sum\limits_{i=1}^{t}x_{v_i}+\sum\limits_{u\in U_0}x_u\le 2x_{u^*}.
\end{align*}
which leads to $\lambda\le 2< 1+\sqrt{m-2}$, a contradiction. Hence, $Y\neq \varnothing$.

Recall that $H^*$ is $3K_2$-free and $H^*[S,T]$ is the maximal bipartite subgraph of $H^*$, so $H^*[Y]$ contains no $P_3$ and $2K_2$. If $H^*[Y]$ contains exactly one edge, denote it by $y_1y_2$, then since $H^*$ is $3K_2$-free and $\delta(H^*)=1$, $Y\setminus\{y_1,y_2\}=\varnothing $,  $d_T(y_1)=d_T(y_2)=1$ and $N_T(y_1)=N_T(y_2)$  which implies that $H^*$ has the form of $D_9$.  Similarly,  by calculation, $\eta_2(D_9)<-4$, a contradiction to (\ref{eql}). Hence, $H^*[Y]$ contains no edge. 

We claim that there are some vertices of $Y$ which have {at least two neighbors} in $T$, Otherwise,  if $d_T(y)=1$ holds for any $y\in Y$, then $H^*$ has the form of $D_{10}$ or $D_{11}$. However,  $\eta_2(D_{10})<-3$ and $\eta_2(D_{11})<-4$, a contradiction. Let $y_1$  be a vertex of $Y$ such that $d_T(y_1)\ge 2$. On the one hand, if {$T=N_T(y_1)$}, then we find a larger dominating complete bipartite subgraph containing $H^*[S,T]$, which contradicts to the maximality of $H^*[S,T]$. So we know that {$N_T(y_1)\subsetneq T$}, which can imply that $t\ge 3$.  On the other hand, if there exists a vertex $y_2\in Y\setminus \{y_1\}$, then since $H^*[S,T]$ is a dominating subgraph, $H^*$ will contains a $3K_2$. Thus, $Y=\{y_1\}$.  Hence, we conclude that  $H^*$ has the form of $D_{12}$. However, one can verify that $\eta_2(D_{12})<-4$, a contradiction to (\ref{eql}).

By the discussion of {Cases} 1 and 2, we conclude that $\delta(H^*)\ge  2$.\QED

\begin{claim}
$G^*= K_3\vee \frac{m-3}{3}K_1$ for $m\ge  33$.
\end{claim}
\noindent \textbf{Proof.} Let $|V(H^*)|=v(H^*)$. Since $\delta(H^*)\ge  2$, we obtain
\begin{align}\label{eqj}
\eta_2(H^*)=\sum\limits_{u\in V(H)}(d_{H^*}(u)-2)\frac{x_u}{x_{u^*}}-e(H^*)\le e(H^*)-2v(H^*),
\end{align}
where the equality occurs if and only if $x_u=x_{u^*}$ for any  $u\in V(H^*)$ with $d_{H^*}(u)\ge 3$.
 We claim that $v(H^*)\ge  6$. Otherwise,  
 if $e(U)=e(H^*)\le 10$, then by Claim \ref{cl-ubo}, 
 \begin{align*}
  {\lambda} \le \frac{1}{2}(e(U)+3)\le \frac{13}{2}<1+\sqrt{m-2},   
 \end{align*}
 a contradiction when $m\ge  33$.
  Therefore,  by Lemma \ref{le-E-G},
  \begin{align}\label{eqk}
   e(H^*)\le 2v(H^*)-3,   
  \end{align}
   and the equality holds if and only if {$H^*=  K_2\vee (v(H^*)-2)K_1$.} Combined it with (\ref{eqj}), we  get  $\eta_2(H^*)\le -3$. Recall that $\eta_2(H^*)\ge  -3$,  so $\eta_2(H^*)=-3$. Note that the equalities in (\ref{eql}), (\ref{eqj}) and (\ref{eqk}) hold,  which implies that

   \begin{itemize}

\item[(d)]  
 $e(W)=0$, 
 $U_0=\varnothing $ and{$H^*= K_2\vee (v(H^*)-2)K_1$};

 \item[(e)] 
 $x_u=x_{u^*}$ for any $u\in U$ with $d_U(u)\ge  3$; 
 
 \item[(f)] 
 $x_w=x_{u^*}$ for any  $w\in W$ with $d_U(w)\ge  1$.  
 \end{itemize}

 We claim that $W=\varnothing $. If the claim holds, then we get $G^*= K_3\vee \frac{m-3}{3}K_1$ immediately.
 Otherwise, if $W\neq \varnothing$, then fix a vertex $w\in W$. Obviously, {$N_{G^*}(w)\subsetneq N_{G^*}(u^*)$}, otherwise, $G^*$ contains a $F_3$. Since $e(W)=0$, we have 
 \[ \lambda x_w={\sum_{u\in N_{G^*}(w)}x_u< \sum_{u\in N_{G^*}(u^*)}x_u}=\lambda x_{u^*}.\]
 It follows that $x_w<x_{u^*}$, which is a contradiction to (f).\QED

\section{Concluding remarks}

In this paper, we have studied the spectral extremal graph problems for graphs of given size with forbidden subgraphs. We determined the spectral extremal graphs for $V_5$ and 
$F_3$. Our result is a unified extension on some recent results on cycles.

In this section, we present some related spectral extremal problems for interested readers. 
Recall that $V_{k+1}= K_1\vee P_{k}$ denotes 
the fan graph on $k+1$ vertices.

\begin{conj}
Let $k\ge 2$ be fixed and $m$ be sufficiently large. 
If $G$ is a $V_{2k+1}$-free or $V_{2k+2}$-free graph with $m$ edges, then 
\[ \lambda(G) \le \frac{k-1 + \sqrt{4m-k^2 +1}}{2}, \]
where the equality holds if and only if $G=K_{k}\vee (\frac{m}{k} - \frac{k-1}{2})K_{1}$. 
\end{conj}

We write $W_{k+1}=K_1\vee C_{k}$ for the wheel graph on $k+1$ vertices.  

\begin{conj}
Let $k\ge 2$ be fixed and $m$ be sufficiently large. 
If $G$ is a $W_{2k+1}$-free  graph with $m$ edges, then 
\[ \lambda(G) \le \frac{k-1 + \sqrt{4m-k^2 +1}}{2}, \]
where the equality holds if and only if $G=K_{k}\vee (\frac{m}{k} - \frac{k-1}{2})K_{1}$. 
\end{conj}

Observe that $W_{2k+2}$ is color-critical and $\chi (W_{2k+2})=4$. 

\begin{conj}
Let $k\ge 2$ be fixed and $m$ be large enough. 
If $G$ is a $W_{2k+2}$-free graph with $m$ edges, then 
\[  \lambda (G) \le \sqrt{4m/3},  \]
with equality if and only if $G$ is a regular complete $3$-partite graph. 
\end{conj}

\section*{Declaration of competing interest}
The authors declare that they have no conflicts of interest to this work.

\section*{Data availability}
No data was used for the research described in the article.

\section*{Acknowledgements} 
The research was supported by the National Natural Science Foundation of China 
(No. 11931002) and the Postdoctoral Fellowship Program of CPSF (No. GZC20233196).

\end{document}